\documentclass[11pt,reqno]{amsart}
\usepackage{lmodern}
\usepackage{microtype}
\usepackage[T1]{fontenc}        
\usepackage[utf8]{inputenc}     
\IfFileExists{stmaryrd.sty}
  {\usepackage{stmaryrd}}
  {\providecommand{\llbracket}{[\![}\providecommand{\rrbracket}{]\!]}}
\usepackage{amssymb}
\usepackage{amscd,mathtools,mathrsfs,amsmath}
\usepackage{amsfonts}
\usepackage{tikz-cd}
\usepackage{tkz-euclide}
\usepackage[all,cmtip]{xy}

\usepackage{xcolor}
\usepackage{hyperref}
\hypersetup{hidelinks,pdfauthor={Takuya Tanaka},pdftitle={Remarks on non-vanishing modulo p of Hecke L-values for imaginary quadratic fields}}
\newcommand\nc{\newcommand}

\theoremstyle{plain}
\newtheorem{theorem}{Theorem}[section]
\newtheorem{prop}[theorem]{Proposition}
\newtheorem{importnota}[theorem]{Important Notation}
\newtheorem{prblm}[theorem]{Problem}
\newtheorem{notation}[theorem]{Notation}
\newtheorem{defin}[theorem]{Definition}
\newtheorem{caution}[theorem]{Caution}
\newtheorem{remark}[theorem]{Remark}
\newtheorem{lemma}[theorem]{Lemma}
\newtheorem{construction}[theorem]{Construction}
\newtheorem{corollary}[theorem]{Corollary}
\newtheorem{example}[theorem]{Example}
\newtheorem{conclusion}[theorem]{Conclusion}
\newtheorem{triviality}[theorem]{Triviality}
\newtheorem{proto}[theorem]{Prototype Quasifibration}
\newtheorem{cauex}[theorem]{Cautionary Example}
\newtheorem{propositiondef}[theorem]{Proposition-Definition}
\newtheorem{subth}{Nuisance}[theorem]
\newtheorem{ssubth}{ }[subth]
\newtheorem{conjecture}[theorem]{Conjecture}
\newtheorem{sidest}[theorem]{Side Story}
\newtheorem{miniexample}[theorem]{Example}

\nc\tri[1]{\begin{triviality}
\label{#1}}
\nc\side[1]{\begin{sidest}
\label{#1}}
\nc\conj[1]{\begin{conjecture}
\label{#1}}
\nc\prodef[1]{\begin{propositiondef}
\label{#1}}
\nc\prt[1]{\begin{proto}
\label{#1}}
\nc\lem[1]{\begin{lemma}
\label{#1}}
\nc\sblm[1]{\begin{sublemma}
\label{#1}}
\nc\pro[1]{\begin{prop}
\label{#1}}
\nc\thm[1]{\begin{theorem}
\label{#1}}
\nc\cor[1]{\begin{corollary}
\label{#1}}
\nc\dfn[1]{\begin{defin}
\label{#1}}
\nc\sthm[1]{\begin{subth}
\label{#1}}
\nc\exm[1]{\begin{example}
\label{#1}}
\nc\miniexm[1]{\begin{miniexample}
\label{#1}}
\nc\plm[1]{\begin{prblm}
\label{#1}}
\nc\rmk[1]{\begin{remark}
\label{#1}}
\nc\subrmk[1]{\begin{subremark}
\label{#1}}
\nc\ntn[1]{\begin{notation}
\label{#1}}
\nc\cau[1]{\begin{caution}
\label{#1}}
\nc\imn[1]{\begin{importnota}
\label{#1}}
\nc\cax[1]{\begin{cauex}
\label{#1}}
\nc\con[1]{\begin{construction}
\label{#1}}
\nc\ssthm[1]{\begin{ssubth}
\label{#1}}
\nc\cnc[1]{\begin{conclusion}
\label{#1}}

\nc\elem{\end{lemma}}
\nc\esblm{\end{sublemma}}
\nc\eside{\end{sidest}}
\nc\econj{\end{conjecture}}
\nc\eprodef{\end{propositiondef}}
\nc\eprt{\end{proto}}
\nc\ethm{\end{theorem}}
\nc\ecor{\end{corollary}}
\nc\edfn{\end{defin}}
\nc\esthm{\end{subth}}
\nc\epro{\end{prop}}
\nc\etri{\end{triviality}}
\nc\eexm{\end{example}}
\nc\eminiexm{\end{miniexample}}
\nc\ermk{\end{remark}}
\nc\subermk{\end{subremark}}
\nc\eplm{\end{prblm}}
\nc\ecau{\end{caution}}
\nc\ecax{\end{cauex}}
\nc\eimn{\end{importnota}}
\nc\entn{\end{notation}}
\nc\econ{\end{construction}}
\nc\ecnc{\end{conclusion}}
\nc\essthm{\end{ssubth}}
\nc\acal{{\mathcal A}}
\nc\affine{{\mathbb A}}
\nc\cc{{\mathbb C}}
\nc\ccc{{\mathcal C}}
\nc\ccal{{\mathcal C}}
\nc\fcal{\mathcal F}
\nc\gcal{{\mathcal G}}
\nc\hcal{{\mathcal H}}
\nc\ical{{\mathcal I}}
\nc\kcal{{\mathcal K}}
\nc\moduli{{\mathcal M}}
\nc\mm{{\bf m}}
\nc\ncal{{\mathcal N}}
\nc\nn{{\bf n}}
\nc\oo{{\mathcal O}}
\nc\pp{{\mathbb P}}
\nc\pcal{{\mathcal P}}
\nc\rr{{\mathbb R}}
\nc\GS{{\mathfrak S}}
\nc\sss{{\mathbb T}}
\nc\xcal{{\mathcal X}}
\nc\icap{{\smallfrown}}

\nc\coker{\mbox{\rm coker}}
\nc\Coker{\mbox{\rm Coker}}
\renewcommand\ker{\mbox{\rm Ker}}
\nc\disjoint{{\amalg}}
\nc\sHom{{\mathcal Hom}}
\nc\cosk{{\rm cosk}}
\nc\tensor{{\otimes}}
\nc\im{\mbox{\rm Im}}
\nc\rank{\mbox{\rm rank}\ }
\nc\codim{\mbox{\rm codim}\ }
\nc\face{\prec}
\nc\Face{\succ}

\newcommand{\C}{\mathbb{C}}

\newcommand{\Z}{\mathbb{Z}}

\renewcommand{\O}{\mathcal{O}}

\def\hoge<#1>{\langle #1 \rangle}

\newcommand{\subsubsubsection}{\@startsection{paragraph}{4}{\z@}
  {1.0\Cvs \@plus.5\Cdp \@minus.2\Cdp}
  {.1\Cvs \@plus.3\Cdp}
  {\reset@font\sffamily\normalsize}
}
\usepackage{chngcntr} 
\title[Remarks on non-vanishing of Hecke $L$-values]{Remarks on non-vanishing modulo $p$ of Hecke $L$-values
  for imaginary quadratic fields}
\author{Takuya Tanaka}
\address{Mathematical Inst. Institute of Science Tokyo\\
 2-12-1, Ookayama, Meguro, Tokyo, 152-8550, JAPAN}
\email{tanaka.t.ch@m.titech.ac.jp}
\AtBeginDocument{\addtolength{\abovedisplayskip}{-1pt}\addtolength{\belowdisplayskip}{-1pt}}
\begin{document}

\begin{abstract}
We extend Lamplugh's nonvanishing theorem to Hecke characters of
arbitrary critical weight over an imaginary quadratic field of class
number one.
\end{abstract}

\maketitle
\tableofcontents

\section{Introduction}

Lamplugh proves a nonvanishing theorem modulo $p$ for critical $L$-values of powers of a Hecke character attached to a CM elliptic curve, along the $\mathbf Z_\ell$-extension of an imaginary quadratic field $K$ unramified outside one of the primes above a split prime $\ell\neq p$ \cite[Theorem~6.9]{Lam15}.
His argument represents the relevant Eisenstein values by logarithmic
derivatives of rational functions and applies his theorem on the
Gamma transforms of elliptic function measures
\cite[Theorem~5.1]{Lam15}. We extend this argument to general critical Hecke characters.
Kundu and Lei \cite{KLv1} studied a generalization of Lamplugh's theorem
concerning $p$-adic valuations of critical Hecke $L$-values in two variables. In their preprint, the
construction in \cite[Definition~6.8 and Appendix~A, (A.2)]{KLv1}
uses $(j,k)$-th Kronecker--Eisenstein functions, which are non-meromorphic for $j < 0$. For $j<0$, the right-hand side of (A.2) is not meromorphic, so the rationality asserted in \cite[Definition~6.8]{KLv1} requires a separate argument. We construct the rational function needed for the
one-variable elliptic measure directly modulo $p$.

Assume throughout that $K$ is an imaginary quadratic field of class
number one and that $E/K$ has complex multiplication by
$\mathcal O_K$. Let $p,\ell\geq5$ be distinct primes which split as
$
 p\mathcal O_K=\mathfrak p\overline{\mathfrak p},
 \ 
 \ell\mathcal O_K=\mathfrak l\overline{\mathfrak l},
$
and assume that $E$ has good reduction at the primes above $p\ell$.
Fix complex and $p$-adic embeddings of $\overline{\mathbf Q}$, the
latter inducing $\mathfrak p$. Choose an invariant differential
$\omega$ generating the module of integral invariant differentials
at $\mathfrak p$.
Define the period lattice and the normalized complex uniformization.

Let $\Lambda=\{\int_\gamma\omega:\gamma\in H_1(E(\C),\Z)\}$.
Let $\pi_{\Lambda}:\C/\Lambda\longrightarrow E(\C)$ be the inverse of the
Abel--Jacobi map $P\mapsto \int_{O}^{P}\omega$. Since $K$ has
class number one, write $\Lambda=\Omega_\infty\mathcal O_K$ with
$\Omega_\infty\in\C^\times$. 
Let $K_{\mathfrak l,\infty}/K$ be the $\mathbf Z_\ell$-extension
unramified outside $\mathfrak l$, and put
$
 \Gamma_{\mathfrak l}
   =\operatorname{Gal}(K_{\mathfrak l,\infty}/K).
$
Fix a critical Hecke character $\lambda$ of conductor
$\mathfrak f_\lambda$ prime to $p\ell$, with the convention
\begin{equation}\label{eq:critical-type}
 \lambda((\alpha))
   =\eta(\alpha)\alpha^{-(a+1)}\overline\alpha^{\,b},
 \qquad a,b\geq0,
 \qquad (\alpha,\mathfrak f_\lambda)=1.
\end{equation}
Set $k=a+1$ and $A_K=\sqrt{|D_K|}/(2\pi)$, where $D_K$ is the
discriminant of $K$. Define the algebraic part of the Hecke $L$--values by
\begin{equation}\label{eq:algebraic-L-value}
 L^{\mathrm{alg}}(\lambda,0)
   =\frac{a!\,L(\lambda,0)}
          {\Omega_\infty^{a+b+1}A_K^b}.
\end{equation}
Its algebraicity follows by expressing the Hecke $L$-value in
terms of Eisenstein--Kronecker numbers and applying
\cite[Corollary~2.10]{BK}. We will prove the following theorem in this paper.

\begin{theorem}\label{thm:main}
Let $\lambda$ have any critical weight as in
\eqref{eq:critical-type}, and let $\varphi$ be a fixed finite-order
Hecke character of conductor dividing $\mathfrak l$.
For all but finitely many finite-order characters $\psi$ of
$\Gamma_{\mathfrak l}$, one has
\begin{equation}\label{eq:main-nonvanishing}
 L^{\mathrm{alg}}(\lambda\varphi\psi,0)
      \in\mathcal O_{\mathbf C_p}^{\times}.
\end{equation}
\end{theorem}
Lamplugh's theorem treats powers of the CM Hecke character in
the range $b=0$ and $0\leq a\leq p-2$, with the hypotheses in
\cite[Theorem~6.9]{Lam15}. In Section~\ref{sec:eisenstein}, we use
integral Poincar\'{e} theta function of Bannai and Kobayashi \cite{BK} to
construct rational functions $\mathcal F_N^{a,b}$ on a $p$-isogenous
curve $q:E'\to E$ with
\[
 \mathcal F_N^{a,b}(Q)
   \equiv E_N^{a,b}(q(Q))\pmod{\mathfrak m_{\mathcal W}}
\]
at nonzero torsion points $Q$ of order prime to $pN$.
Section~\ref{sec:measures} uses translates of these functions to
construct an elliptic function measure and proves that its restriction
to the units has nonzero reduction. Then, Lamplugh's Gamma-transform theorem proves Theorem~\ref{thm:main}.
\section{Congruences for Eisenstein--Kronecker numbers}
\label{sec:eisenstein}

Retain the notation of the introduction, and let $\mathcal E_{/\O_{K,\mathfrak p}}$ be
the N\'{e}ron model of $E$ at $\mathfrak p$.
Use the identification of the CM action with $\mathcal O_K$ for
which $[\alpha]^*\omega=\alpha\omega$. Then, the finite group scheme $\mathcal E[\mathfrak p]$
is connected. Put $\mathcal W_0=W(\overline{\mathbf F}_p)$ and
$\mathcal W=\mathcal W_0[\eta,\mu_{\ell-1}]$, adjoining the values
of $\eta$ and the roots of unity of order dividing $\ell-1$.
Let $J=\operatorname{Frac}(\mathcal W)$,
$\mathbf k=\overline{\mathbf F}_p$, and let $\mathfrak m_{\mathcal W}$
be the maximal ideal of $\mathcal W$. For $u\in\mathbf k^\times$,
write $[u]_{\mathrm{Teich}} \in\mathcal W_0^\times$ for its
Teichm\"uller lift. 

\subsection{Eisenstein measure}

Put $A_\Lambda=\frac{\operatorname{Area}(\mathbf C/\Lambda)}{\pi}=|\Omega_\infty|^2A_K$, and 
 \[
 \langle ,\rangle_\Lambda:\C\times\C\rightarrow \C^\times:\langle z,w\rangle_\Lambda
 =\exp\!\left(\frac{z\overline w-w\overline z}{A_\Lambda}\right).
\]
We use the following complex function on $(z,w,u)\in \C\times \C\times \C$ which defined in \cite[\S1.1]{BK}:
\begin{align*}
 K_m^*(z,w,u;\Lambda)
 &=
 \sum_{\substack{\gamma\in\Lambda\\ \gamma\ne-z}}
 \frac{\overline{z+\gamma}^{\,m}}{|z+\gamma|^{2u}}
 \langle\gamma,w\rangle_\Lambda,
 \qquad \operatorname{Re}(u)>\frac m2+1,\\
 e^*_{r,s}(z,w;\Lambda)
 &=K^*_{r+s}(z,w,s;\Lambda)
 \qquad (r\geq0,\ s\geq1).
\end{align*}
Values outside the region of absolute convergence are defined by
analytic continuation.
Fix $N\geq2$ with $(N,p)=1$, and put
$f=N^2\mathbf1_{\{O\}}-\mathbf1_{E[N]}$ and $c_t=f(-t)$.
Then, we have $\sum_{t\in E[N]}c_t=0$.

\begin{defin}
For $a,b\geq0$ and torsion point $x\neq 0$ with its order prime to $pN$, define
\begin{equation}\label{eq:BK-EK-definition}
 E_N^{a,b}(x)
 =\frac{a!}{A_\Lambda^b}
    \sum_{t\in E[N]}c_t e^*_{b,a+1}(x+t,0;\Lambda).
\end{equation}
These values are algebraic by \cite[Corollary~2.10]{BK}; we view
them $p$-adically through the fixed embedding.
\end{defin}
Let $\widehat{\mathcal{E}}$, and $\widehat{\mathbf G}_m$ be the formalizations along their unit section, and choose a formal group isomorphism over $\mathcal W_0$ 
\[
 \xi:\widehat{\mathbf G}_m\xrightarrow{\sim}\widehat{\mathcal E},
 \qquad
 \log_E(\xi(T))=\Omega_p\log(1+T),
 \qquad \Omega_p\in\mathcal W_0^\times,
\]
where $\log_E$ is the formal logarithm associated with $\omega$;
see \cite[\S3.1]{BK}.
Let $\Theta_{z_0,w_0}(z,w;\Lambda)$ be the translated theta function
of \cite[\S1.3]{BK} and write $\widehat\Theta_{z_0,w_0}$ for its
expansion after substituting $\log_E(s),\log_E(t)$ for $z,w$. Here, $s,t$ are a copy of two fixed local parameters of $\mathcal{E}$ at $O$ such that $\log_E(t)=t+O(t^2)$. Put
\[
 H_{N,x}(S,T)=\sum_{\tau\in E[N]}c_\tau
      \widehat\Theta_{x+\tau,0}(\xi(S),\xi(T)).
\]
Here, $x\neq 0$ is the fixed torsion point with its order coprime to $pN$.
\begin{prop}\label{prop:smoothed-measure}
The series $H_{N,x}$ belongs to $\mathcal W\llbracket S,T\rrbracket$ and is the
Amice transform of a unique $\mathcal W$-valued measure
$\mu_{N,x}$ on $\mathbf Z_p^2$. Its moments satisfy
\begin{equation}\label{eq:smoothed-moments}
 E_N^{a,b}(x)
  =\frac{(-1)^{a+b}}{\Omega_p^{a+b}}
       \int_{\mathbf Z_p^2}u^av^b\,d\mu_{N,x}(u,v)
  \in\mathcal W.
\end{equation}
\end{prop}

\begin{proof}
 By \cite[Corollary~2.17]{BK},
$\widehat\Theta_{x+\tau,0}(s,t)-t^{-1}$ is integral for
$\tau\in E[N]$. The relation $\sum_\tau c_\tau=0$ cancels the common
pole, so $H_{N,x}$ is integral. By
\cite[Theorem~1.17]{BK}, we have
\[
 \sum_t c_t\Theta_{x+t,0}(z,w)
  =\sum_{a,b\geq0}\frac{(-1)^{a+b}}{b!A_\Lambda^b}
       \left(\sum_t c_t e^*_{b,a+1}(x+t,0;\Lambda)\right)z^aw^b.
\]
Apply the Amice correspondence \cite[Proposition~3.1]{BK},
followed by $((1+S)\partial_S)^a((1+T)\partial_T)^b$ at the
origin. This gives \eqref{eq:smoothed-moments}.
\end{proof}

\subsection{A rational function on the isogenous curve}

Let $q:E'\longrightarrow E$ be the degree-$p$ isogeny dual to
$E\longrightarrow E/E[\mathfrak p]$. Thus $E'=E/E[\mathfrak p]$
and $q$ is \'{e}tale on the integral models. Put $\omega'=q^*\omega$.
The corresponding complex uniformization is
\[
 E'(\mathbf C)=\mathbf C/\Lambda',
 \qquad \Lambda'=\overline{\mathfrak p}\Lambda,
\]
and $q$ is induced by the identity on $\mathbf C$.
Write $C=\ker(q)=\Lambda/\Lambda'$ and
$H=E[\mathfrak p]=\mathfrak p^{-1}\Lambda/\Lambda$.
The pairing
\[
 C\times H\longrightarrow\mu_p,
 \qquad (s,u)\longmapsto\langle s,u\rangle_\Lambda
\]
is well-defined and perfect.

Fix a primitive $p$-th root of unity $\zeta_p$ and use the following convention for each $0\leq r<p$.
\[
 \xi'=\widehat q^{-1}\circ\xi,
 \qquad u_r=\xi(\zeta_p^r-1),
 \qquad u'_r=\xi'(\zeta_p^r-1)\qquad (0\leq r<p).
\]
For $0\leq j<p$, let $s_j\in C$ be the element determined by
\begin{equation}\label{eq:dual-residue-labels}
 \langle s_j,u_r\rangle_\Lambda=\zeta_p^{jr}
 \qquad(0\leq r<p).
\end{equation}
Each prime-to-$p$ torsion point has a unique prime-to-$p$ lift
through $q$; write $t'$ and $x'$ for the lifts of $t$ and $x$.

Choose a minimal integral Weierstrass equation
\begin{equation}\label{eq:Weierstrass-model}
 E':Y^2=4X^3-g_2X-g_3,
 \qquad \omega'=\frac{dX}{Y},
 \qquad g_2,g_3\in\mathcal O_J,
\end{equation}
and denote its reduction by $\widetilde E'$.
Normalize $v_J$ by $v_J(\varpi_J)=1$ where $\varpi_J$ is the uniformizer of $J$.

\begin{defin}\label{def:Lamplugh-val}
By \cite[Lemma~3.2, (3.1), and the following remark]{Lam15},
every nonzero $h\in J(E')$ can be written as
\[
 h=\varpi_J^r\frac{F(X,Y)}{G(X,Y)},
 \qquad r\in\mathbf Z,\quad F,G\in\mathcal O_J[X,Y],
\]
where $\overline F$ and $\overline G$ are nonzero functions on
$\widetilde E'$. For every prime $\ell\ne p$, both $F(P)$ and
$G(P)$ are units for all but finitely many
$P\in E'[\ell^\infty]$.
We define $v(h)=r$ and $v(0)=+\infty$.
This is a valuation on $J(E')$ extending $v_J$.
\end{defin}

\begin{lemma}\label{lem:R_f}
The expression
\begin{equation}\label{eq:rational-base-function}
 R_f(Q)=\sum_{t\in E[N]}c_t e^*_{0,1}(Q+t',0;\Lambda')
\end{equation}
defines a rational function on $E'$ over a finite unramified extension of $\mathbf Q_p$.
Its poles are simple, at each $-t'\in E'[N]$, with residues $c_t$ relative to
$\omega'$, and $v(R_f)=0$.
\end{lemma}

\begin{proof}
With $X=\wp(z)$ and $Y=\wp'(z)$, recall that
\[
 \zeta(z)=e^*_{0,1}(z,0;\Lambda')
           +e_2^*(\Lambda')z+\frac{\overline z}{A_{\Lambda'}},
 \qquad e_2^*(\Lambda')=e^*_{0,2}(0,0;\Lambda');
\]
see \cite[Example~1.9, (15), and Theorem~1.17]{BK}. Here, $\zeta(z)$ is the Weierstrass $\zeta$-function.
The correction term is additive, so the Weierstrass addition
formula \cite[(23.10.2)]{DLMF}
\[
 \zeta(z+w)=\zeta(z)+\zeta(w)
       +\frac{\wp'(z)-\wp'(w)}{2(\wp(z)-\wp(w))}
\]
gives the same addition formula for $e^*_{0,1}$.
Since $\sum_t c_t=0$, it follows that
\begin{equation}\label{eq:base-addition}
 R_f(Q)=\sum_{\substack{t\in E[N]\\t\ne O}}c_t
   \left(e^*_{0,1}(t',0;\Lambda')
       +\frac{Y(Q)-Y(t')}{2(X(Q)-X(t'))}\right).
\end{equation}
The points $O\ne t'\in E'[N]$ have order prime to $p$, so they lie in
finite unramified extensions of the local field of the model.
By \cite[Corollary~2.10]{BK} the constants lie in the field
obtained by adjoining $E'[4N^2]$ to the field of definition
of the model. This is still unramified at $p$. Moreover, by
\cite[Corollary~2.17]{BK}, $e^{*}_{0,1}(t',0;\Lambda')$ for each $t'\in E'[N]$ are $p$-adically integral.
The fractions $\frac{Y(Q)-Y(t')}{2(X(Q)-X(t'))}$ have valuation zero, hence
$v(R_f)\geq0$. The simple poles at distinct prime-to-$p$ torsion
points remain distinct after reduction, and the residue at
$-t'$ is $c_t=-1$ for $O\ne t\in E[N]$. Thus, we have $v(R_f)=0$.
\end{proof}

\begin{prop}\label{prop:mass-formula}
For $0\leq i,j<p$, put
$U_{i,j}=(i+p\mathbf Z_p)\times(j+p\mathbf Z_p)$.
For every torsion point $x\neq 0$ with its order prime to $pN$,
\begin{equation}\label{eq:mass-formula}
 \mu_{N,x}(U_{i,j})
   =\frac1p\sum_{r=0}^{p-1}
        \zeta_p^{-ir}R_f(x'+u'_r+s_j).
\end{equation}
In particular, the right-hand side is integral.
\end{prop}

\begin{proof}
First, for $w\in\mathfrak p^{-1}\Lambda/\Lambda$, decomposition
into cosets of $\Lambda'$ gives
\begin{equation}\label{eq:lattice-decomposition}
 e^*_{0,1}(z,w;\Lambda)
  =\sum_{s\in\Lambda/\Lambda'}
       \langle s,w\rangle_\Lambda
       e^*_{0,1}(z+s,0;\Lambda').
\end{equation}
This follows first for the convergent series $K_1^*$ and then
by analytic continuation.

We next claim that, for every $0\leq r,v<p$, 
\begin{equation}\label{eq:theta-torsion-values}
 H_{N,x}(\zeta_p^r-1,\zeta_p^v-1)
    =\sum_{j=0}^{p-1}\zeta_p^{jv}R_f(x'+u'_r+s_j).
\end{equation}
Let $n$ be a positive integer which is prime to $p$ and annihilate $x$, $E[N]$. Choose
$\epsilon\in\mathcal O_K$ with
\[
 \epsilon\equiv1\pmod{n\mathcal O_K},
 \qquad \epsilon\equiv1\pmod{\mathfrak p},
 \qquad \epsilon\equiv0\pmod{\overline{\mathfrak p}}.
\]
For each $\tau\in E[N]$, choose $z_\tau\in n^{-1}\Lambda'$
representing $x'+\tau'$. The translation formula
\cite[Corollary~2.21]{BK} leads
\[
 \widehat\Theta_{z_\tau,0}(s\oplus u_r,t\oplus u_v)
   =\langle\epsilon u_v,(\epsilon-1)z_\tau\rangle_\Lambda
      \widehat\Theta_{z_\tau+\epsilon u_r,\,\epsilon u_v}(s,t).
\]
Its pairing factor is $1$. Actually,
$(\epsilon-1)z_\tau\in\Lambda'$, which annihilates
$\epsilon u_v\in\mathfrak p^{-1}\Lambda/\Lambda$.
Sum with coefficients $c_\tau$ and evaluate at the origin, canceling
the common pole of the second variable when $v=0$. This gives
\[
 H_{N,x}(\zeta_p^r-1,\zeta_p^v-1)
    =\sum_{\tau\in E[N]}c_\tau
       e^*_{0,1}(z_\tau+\epsilon u_r,\epsilon u_v;\Lambda).
\]
Now $\epsilon u_r\bmod\Lambda'=u'_r$ and
$\langle s_j,\epsilon u_v\rangle_\Lambda=\zeta_p^{jv}$.
Applying \eqref{eq:lattice-decomposition} proves
\eqref{eq:theta-torsion-values}.

Finally, the identity
$\mathbf1_{i+p\mathbf Z_p}(t)
 =p^{-1}\sum_{r=0}^{p-1}\zeta_p^{r(t-i)}$ gives
\[
 \mu_{N,x}(U_{i,j})
   =\frac1{p^2}\sum_{r,v=0}^{p-1}
      \zeta_p^{-ir-jv}H_{N,x}(\zeta_p^r-1,\zeta_p^v-1).
\]
Insert \eqref{eq:theta-torsion-values} and sum over $v$ to obtain
\eqref{eq:mass-formula}.
\end{proof}

\begin{corollary}\label{cor:congru-Eis-rational}
For $a,b\geq0$ and torsion point $x\neq 0$ with its order prime to $pN$,
\begin{equation}\label{eq:finite-EK-congruence}
 E_N^{a,b}(x)
  \equiv\frac{(-1)^{a+b}}{\Omega_p^{a+b}}
       \sum_{i,j=0}^{p-1}i^aj^b
       \left(\frac1p\sum_{r=0}^{p-1}
          \zeta_p^{-ir}R_f(x'+u'_r+s_j)\right)
       \pmod{\mathfrak m_{\mathcal W}}.
\end{equation}
The summation is formed before reduction.
\end{corollary}

\begin{proof}
On $U_{i,j}$, one has $u^av^b\equiv i^aj^b\pmod p$.
Apply \eqref{eq:smoothed-moments} and \eqref{eq:mass-formula}.
\end{proof}
Let $D$ be the invariant derivation dual\footnote{So, for each rational function $h$, $Dh$ is a rational function defined by $dh=Dh\cdot  \omega^{'}$.} to $\omega'$.
In the Weierstrass coordinates given in \eqref{eq:Weierstrass-model}, $D$ is explicitly denoted by
\[
 D=Y\frac{\partial}{\partial X}
      +\left(6X^2-\frac{g_2}{2}\right)\frac{\partial}{\partial Y}.
\]
Thus, $D$ sends integral coordinate polynomials to integral
polynomials. For $h=\varpi_J^r F/G$ as in
Definition~\ref{def:Lamplugh-val},
\[
 Dh=\varpi_J^r\frac{(DF)G-F(DG)}{G^2},
 \qquad v(Dh)\geq v(h).
\]
Put $\kappa_0=\overline{\Omega_p}\in\mathbf k^\times$ and let $\widetilde\kappa_0=[\kappa_0]_{\mathrm{Teich}}$ be its Teichmuller lift.

\begin{lemma}\label{lem:ell-func}
For $h\in J(E')$ with $v(h)\geq0$ and $a\geq0$, put
\[
 \mathcal C_a h(Q)
   =\frac1p\sum_{i,r=0}^{p-1}
           i^a\zeta_p^{-ir}h(Q+u'_r).
\]
Then, $v(\mathcal C_a h)\geq0$ and we have
 $\overline{\mathcal C_a h}=\delta^a\overline h$, where $\delta$ be the derivative defined by $\delta=\kappa_0\overline D$.
Moreover, $\delta^p=\delta$ on $\mathbf k(\widetilde E')$.
\end{lemma}

\begin{proof}
Initially $\mathcal C_a h$ is defined over $J(\zeta_p)$. If a
Galois automorphism over $J$ sends $\zeta_p$ to $\zeta_p^d$, it
sends $u'_r$ to $u'_{dr}$, since $\xi'$ is defined over $J$.
Reindexing the sum by $r\mapsto dr\bmod p$ shows that
$\mathcal C_a h\in J(E')$.

Both $\mathcal C_a$ and $D$ commute with translations. For
$h=\varpi_J^r F/G\ne0$, choose $P$ in the
prime to $p$ torsion of $E'$ such that reduction
$\overline P$ avoids the poles of the coordinates and the zeros
of $\overline G$. There are infinitely many such points with distinct
reductions, since this torsion is \'{e}tale. Translate $h$ by $P$
so that $h$ is regular at the origin of the integral model.
Then $h(\xi'(T))\in\mathcal W\llbracket T\rrbracket$ and by the
Amice correspondence, there is a measure $\mu_h$ with
\[
 h(\xi'(T))=\int_{\mathbf Z_p}(1+T)^t\,d\mu_h(t),
 \qquad
 (\Omega_pD)^a h(O)=\int_{\mathbf Z_p}t^a\,d\mu_h(t).
\]
The identity
$\mathbf1_{i+p\mathbf Z_p}(t)
 =p^{-1}\sum_{r=0}^{p-1}\zeta_p^{r(t-i)}$ gives
\[
 \mathcal C_a h(O)
    =\sum_{i=0}^{p-1}i^a\mu_h(i+p\mathbf Z_p)
    \equiv(\Omega_pD)^a h(O)\pmod{p\mathcal W}.
\]
Likewise, $t^p\equiv t\pmod p$ gives
$((\Omega_pD)^p-\Omega_pD)h(O)\in p\mathcal W$.
Undoing the translation gives these statements at infinitely
many $P$ with distinct reductions. Lamplugh's valuation and
the identity principle for rational functions on
$\widetilde E'$ prove the assertions.
\end{proof}

\begin{theorem}\label{thm:rational-representative}
Put
\begin{equation}\label{eq:rational-representative}
 h_b(Q)=\sum_{j=0}^{p-1}j^b R_f(Q+s_j),
 \qquad
 \mathcal F_N^{a,b}(Q)
   =(-1)^{a+b}\widetilde\kappa_0^{-b}D^a h_b(Q).
\end{equation}
Then $\mathcal F_N^{a,b}\in J(E')$ and
$v(\mathcal F_N^{a,b})=0$.
For every prime to $pN$ torsion point $Q\neq 0$ of $E^{'}$, both sides below are defined and integral, and
\begin{equation}\label{eq:rational-interpolation}
 \mathcal F_N^{a,b}(Q)
     \equiv E_N^{a,b}(q(Q))\pmod{\mathfrak m_{\mathcal W}}.
\end{equation}
\end{theorem}

\begin{proof}
Lemma~\ref{lem:R_f} and $v(Dh)\geq v(h)$ show that
$\mathcal F_N^{a,b}$ is integral.
For $Q$ as in the statement, $h_b$ is regular on $Q$. Corollary~\ref{cor:congru-Eis-rational}
and the translated Amice argument of Lemma~\ref{lem:ell-func} give
\[
 \begin{aligned}
 E_N^{a,b}(q(Q))
   &\equiv(-1)^{a+b}\Omega_p^{-a-b}\mathcal C_a h_b(Q)\\
   &\equiv(-1)^{a+b}\Omega_p^{-b}D^a h_b(Q)
    \equiv\mathcal F_N^{a,b}(Q)
        \pmod{\mathfrak m_{\mathcal W}}.
 \end{aligned}
\]

To prove equality of valuations, choose $O\ne t_0\in E[N]$.
The reductions of the points $-t'-s_j$ for each $t'\in E'[N]$, $s_j\in E'[\overline{\mathfrak{p}}]$ are pairwise distinct,
so $\overline{h_b}$ has a simple pole at
$-\overline{t_0'}-\overline{s_1}$ with nonzero residue.
For $a>0$, choose $1\leq a_0\leq p-1$ with
$a_0\equiv a\pmod{p-1}$. Since $\delta^p=\delta$,
one has $\delta^a=\delta^{a_0}$.
In a local parameter $z$ at this pole, write
$\overline{h_b}=cz^{-1}+O(1)$ and $\delta z=u+O(z)$,
where $c,u\in\mathbf k^\times$. Then, we have
\[
 \delta^a\overline{h_b}
   =(-1)^{a_0}a_0!\,c\,u^{a_0}z^{-a_0-1}
       +O(z^{-a_0})\ne0.
\]
For $a=0$, the simple pole itself proves nonvanishing.
Since $\overline{\mathcal F_N^{a,b}}$ is a nonzero scalar multiple
of $\delta^a\overline{h_b}$, we obtain
$v(\mathcal F_N^{a,b})=0$. Here, 
$\mathcal F_N^{a,b}$ is defined over a finite unramified extension of
$\mathbf Q_p$. Actually, this follows from Lemma~\ref{lem:R_f} and the points $s_j\in E'[\overline{\mathfrak p}]$ are defined over a finite unramified extension.
\end{proof}

\section{Elliptic function measures and non-vanishing}\label{sec:measures}
Retain the notation of the preceding sections and put
$w_K=\#\mathcal O_K^\times$. Choose a nonzero proper ideal
$\mathfrak g=(\gamma_\lambda)$ divisible by $\mathfrak f_\lambda$ and prime
to $p\ell$, and choose $N\geq2$ prime to $p\ell\mathfrak g$. So, if $\mathfrak{f}_\lambda\neq (1)$, we may choose $\mathfrak{g}=\mathfrak{f}_\lambda$.
Regard $\eta$ in \eqref{eq:critical-type} as a character modulo
$\mathfrak g$ and put $\eta_0(t)=[\overline{\eta(t)}]_{\mathrm{Teich}}$ be the Teichmuller lift of the reduction.
For $u\in(\mathcal O_K/\mathfrak g)^\times$, let
$V_u'$ be the unique prime-to-$p$ lift through $q$ of
$V_u=\pi_{\Lambda}(\Omega_\infty u/\gamma_\lambda\pmod\Lambda)\in E[\mathfrak g]$, and define the rational functions over $J$ as follows;
\begin{equation}\label{eq:lm-theta}
 \begin{aligned}
 \Theta(P)
   &=\frac1{w_K}
      \sum_{u\in(\mathcal O_K/\mathfrak g)^\times}
         \eta_0(u)\mathcal F_N^{a,b}(P+V_u'),\\
 \Theta^\circ(P)
   &=\frac1{w_K}\sum_{\varepsilon\in\mathcal O_K^\times}
                           \Theta([\varepsilon]P),\\
 R(P)
   &=\Theta^\circ(P)
      -\frac1\ell\sum_{T\in E'[\mathfrak l]}
                           \Theta^\circ(P+T).
 \end{aligned}
\end{equation}
Write $\mathfrak m$ for the maximal
ideal of $\mathcal W$.

\begin{lemma}\label{lem:primitive-kernel}
The functions in \eqref{eq:lm-theta} are integral at every point of
$E'[\mathfrak l^\infty]$. Moreover,
\[
 \overline{\Theta^\circ}=\overline\Theta,
 \qquad v(\Theta^\circ)=v(R)=0.
\]
Both $\Theta^\circ$ and $R$ are invariant under
$\mathcal O_K^\times$.
\end{lemma}

\begin{proof}
For $P\in E'[\mathfrak l^\infty]$, the point $q(P)+V_t$ is
nonzero of order prime to $pN$. Thus,
Theorem~\ref{thm:rational-representative} applies to every summand,
including when $P=O$. The identities
\[
 E_N^{a,b}(\varepsilon x)
   =\varepsilon^{-k}\overline\varepsilon^{\,b}E_N^{a,b}(x),
 \qquad
 \eta(\varepsilon)=\varepsilon^k\overline\varepsilon^{-b}
\]
show that $\overline\Theta$ is invariant under the units, proving
$\overline{\Theta^\circ}=\overline\Theta$. Here, we use the notation $k=a+1$.

The possible poles of $\overline\Theta$ are
$-t'-s_j-V_u'$, where $t\in E[N]$ and $s_j\in E^{'}[\overline{\mathfrak p}]$, $u\in (\O_K/\mathfrak g)^\times$. The $N$-, $p$-, and $\mathfrak g$-torsion subgroups are
independent, and reduction is injective on their sum here. The pole of
$\overline{\mathcal F_N^{a,b}}$ at $-t_0'-s_1$, with fixed $O\neq t_0\in E[N]$,
therefore survives in each translate with a unit coefficient.
Hence, $v(\Theta^\circ)=0$.
The poles of $\overline{\Theta^\circ}$ are torsion points of order prime to $\ell$, so its distinct translates by $E'[\mathfrak l]$ have disjoint pole sets.
Each term indexed by $T\ne O$ therefore contributes a pole to $\overline R$ with coefficient $-1/\ell$, giving $v(R)=0$.
The invariance of $R$ under the action of $\O_K^\times$ follows from the definition of \eqref{eq:lm-theta}.
\end{proof}

\subsection{The measure and its interpolation formula}
Choose an isomorphism
$\mathbf e_\ell:\mathbf Q_\ell/\mathbf Z_\ell
\longrightarrow\mu_{\ell^\infty}$.
Identify $K_{\mathfrak l}/\O_{K,\mathfrak l}
   =\bigcup_{m\geq0}\mathfrak l^{-m}/\mathcal O_K$.
For each $z\in K_{\mathfrak l}/\O_{K,\mathfrak l}$, let $P_z\in E'[\mathfrak l^\infty]$ be a torsion point characterized by
$q(P_z)=\pi_{\Lambda}(\Omega_\infty z\pmod\Lambda)$.
Define a bounded integral measure $\alpha$
on $\O_{K,\mathfrak l}$ by
\begin{equation}\label{eq:lm-measure}
 \alpha(u+\ell^m\O_{K,\mathfrak l})
   =\frac1{\ell^m}
      \sum_{z\in\mathfrak l^{-m}/\mathcal O_K}
         \mathbf{e}_\ell(-uz)\Theta^\circ(P_z).
\end{equation}
We can check the distribution properties of
\eqref{eq:lm-measure} by a simple computation. Fourier transform $\widehat{\alpha}:K_{\mathfrak{l}}/\O_{K,\mathfrak{l}}\rightarrow J;x\mapsto \int_{\O_{K,\mathfrak{l}}}\mathbf{e}_\ell(xt)d\alpha(t)$ of the measure $\alpha$ satisfy $\widehat\alpha(z)=\Theta^\circ(P_z)$.
The boundedness of the measure $\alpha$ follows from Lemma~\ref{lem:primitive-kernel} and
$\ell\ne p$. Thus, $\alpha$ is an elliptic function measure for $E'/J$ in the sense of \cite[\S3.2, p.~616]{Lam15}. 
Its restriction $\mu=\alpha|_{\O_{K,\mathfrak l}^{\times}}$ has Fourier transform $R(P_z)$ and is invariant under the pull-back of $\mathcal O_K^\times$.

Let $\chi$ be a non-trivial finite-order Hecke character of  conductor
$\mathfrak f_\chi=(\nu_\chi)=\mathfrak l^m$, and put
$\rho_\chi(u)=\chi((u))$ on $(\mathcal O_K/\mathfrak f_\chi)^\times$.
Define
\begin{equation}\label{eq:lm-gauss}
 \tau(\rho_\chi^{-1})
   =\frac1{\ell^m}
      \sum_{u\in(\mathcal O_K/\mathfrak f_\chi)^\times}
        \rho_\chi(u)^{-1}\mathbf e_\ell(-u/\nu_\chi).
\end{equation}
The additive character in this sum is primitive. Note that $\tau(\rho_\chi)\tau(\rho_\chi^{-1})=\rho_\chi(-1)\ell^{-m}$, so
$\tau(\rho_\chi^{-1})$ is a $p$-adic unit. Write $L_{\mathfrak g}^{\mathrm{alg}}$ for
the algebraic part of the Hecke $L$-values \eqref{eq:algebraic-L-value} with the Euler factors
at $\mathfrak g$ omitted, denote $k=a+1$ like before and set
\begin{equation}\label{eq:lm-factors}
 \begin{aligned}
 C_\lambda(\chi)
   &=\tau(\rho_\chi^{-1})
      \frac{\gamma_\lambda^k}{\overline\gamma_\lambda^{\,b}}
      \lambda(\mathfrak f_\chi)^{-1}\chi(\mathfrak g)^{-1},\\
 B_N(\lambda\chi)
   &=N^2-(\lambda\chi)((N))^{-1}.
 \end{aligned}
\end{equation}
In particular, $C_\lambda(\chi)$ is a $p$-adic unit.

\begin{prop}\label{prop:lm-interpolation}
The measure $\mu$ satisfies
\begin{equation}\label{eq:lm-interpolation}
 \Gamma_\mu(\rho_\chi^{-1})
  :=\int_{\O_{K,\mathfrak{l}}^\times}\rho_\chi(u)^{-1}\,d\mu(u)
  \equiv C_\lambda(\chi)B_N(\lambda\chi)
       L_{\mathfrak g}^{\mathrm{alg}}(\lambda\chi,0)
       \pmod{\mathfrak m}.
\end{equation}
\end{prop}

\begin{proof}
The distribution relation gives
\[
 E_N^{a,b}(z)
   =\frac{a!}{A_\Lambda^b}
      \left(N^2e_{b,k}^{*}(z,0;\Lambda)
            -N^{k-b}e_{b,k}^{*}(Nz,0;\Lambda)\right).
\]
Set $d=\gamma_\lambda\nu_\chi\in \O_K$, where $\mathfrak g=(\gamma_\lambda)$, $\mathfrak f_\chi=(\nu_\chi)$. The same computation as in
\cite[Theorem~6.2]{Lam15} gives
\begin{equation}\label{eq:lm-ray-sum}
 \begin{aligned}
 &\frac1{w_K}\sum_{r\in(\mathcal O_K/(d))^\times}
       \eta(r)\rho_\chi(r)
       E_N^{a,b}(\Omega_\infty r/d)\\
 &\qquad
   =\frac{d^k}{\overline d^{\,b}}
      B_N(\lambda\chi)L_{\mathfrak g}^{\mathrm{alg}}(\lambda\chi,0).
 \end{aligned}
\end{equation}
Theorem~\ref{thm:rational-representative} and
Lemma~\ref{lem:primitive-kernel} turn \eqref{eq:lm-ray-sum} into
\begin{equation}\label{eq:lm-torsion-sum}
 \begin{aligned}
 &\sum_{u\in(\mathcal O_K/\mathfrak f_\chi)^\times}
       \rho_\chi(u)\Theta^\circ([u]P_{1/\nu_\chi})\\
 &\qquad\equiv
   \frac{\gamma_\lambda^k}{\overline\gamma_\lambda^{\,b}}
      \lambda(\mathfrak f_\chi)^{-1}\chi(\mathfrak g)^{-1}
      B_N(\lambda\chi)L_{\mathfrak g}^{\mathrm{alg}}(\lambda\chi,0)
      \pmod{\mathfrak m}.
 \end{aligned}
\end{equation}
We will prove this shortly. For $t\in (\O_K/\mathfrak g)^\times$, $u\in (\O_K/\mathfrak f_\chi)^\times$, put $Q_{t,u}:=[u]P_{1/\nu_\chi}+V_t^{'}$. By definition, we have $q(Q_{t,u})=\pi_{\Lambda}\!\left(\Omega_\infty(u/\nu_\chi+t/\gamma_\lambda)\bmod\Lambda\right)\allowbreak=\pi_{\Lambda}\!\left(\Omega_\infty\frac{u\gamma_\lambda+t\nu_\chi}{d}\bmod\Lambda\right)$. By the Chinese Remainder Theorem, the set $\{u\gamma_\lambda+t\nu_\chi\in (\O_K/d)^\times\ |\ u\in (\O_K/\mathfrak f_\chi)^\times, t\in (\O_K/\mathfrak g)^\times\}$ is the complete representatives of $(\O_K/d)^\times$. Therefore, the left-hand side of \eqref{eq:lm-torsion-sum} is computed as follows;
\begin{equation}
\begin{aligned}
&\equiv
  \frac{1}{w_K}
  \sum_{u\in(\mathcal O_K/\mathfrak f_\chi)^\times}
    \rho_\chi(u)
  \sum_{t\in(\mathcal O_K/\mathfrak g)^\times}
    \eta(t)\,
    E_N^{a,b}\bigl(q(Q_{t,u})\bigr)
  \pmod{\mathfrak m}
\\
&=
  \rho_\chi(\gamma_\lambda)^{-1}
  \eta(\nu_\chi)^{-1}
  \Biggl(
    \frac{1}{w_K}
    \sum_{r\in(\mathcal O_K/(d))^\times}
      \rho_\chi(r)\eta(r)\,
      E_N^{a,b}\!\left(\Omega_\infty\,\frac{r}{d}\right)
  \Biggr)
\\
&=
  \rho_\chi(\gamma_\lambda)^{-1}
  \eta(\nu_\chi)^{-1}\,
  B_N(\lambda\chi)\,
  \frac{d^k}{\overline d^{\,b}}\,
  L_{\mathfrak g}^{\mathrm{alg}}(\lambda\chi,0).
\end{aligned}
\end{equation}
By using $\eta(\nu_\chi)^{-1}d^k/\overline{d}^b=\lambda(\mathfrak f_\chi)^{-1}\gamma_\lambda^k/\overline{\gamma_\lambda}^b$, we conclude our equation.
Finally, extend $\rho_\chi$ by zero on the nonunits modulo $\mathfrak l^m$. By primitivity and \eqref{eq:lm-gauss},
\[
 \sum_{a\in(\O_K/\mathfrak l^m)^\times}\rho_\chi(a)^{-1}
       \mathbf e_\ell(-av/\nu_\chi)
   =\ell^m\tau(\rho_\chi^{-1})\rho_\chi(v)
 \qquad(v\in\O_K/\mathfrak l^m).
\]
Substituting this identity in \eqref{eq:lm-measure} gives \eqref{eq:lm-interpolation} by \eqref{eq:lm-torsion-sum} (cf.~\cite[Lemma~3.1]{Lam15}).
\end{proof}

\subsection{Proof of Theorem~\ref{thm:main}}
Let $\rho_\varphi:(\O_K/\mathfrak l)^\times\to\mathbf C^\times$ be given by $\rho_\varphi(a)=\varphi((a))$. Replacing $\rho_\varphi$ by $[\overline{\rho_\varphi}]_{\mathrm{Teich}}:(\O_K/\mathfrak l)^\times\rightarrow \overline{\mathbf F_p}^\times\rightarrow \mathcal{W}^\times$, which is the Teichm\"uller lift of the reduction, leaves $\Gamma_\mu(\rho_{\varphi\psi}^{-1})$ unchanged modulo $\mathfrak m$ for every $\psi$, since the two characters have the same reduction. Thus, for the unit assertion in \eqref{eq:gamma-nonvanishing}, we may assume that $\rho_\varphi$ has order prime to $p$.
For the measure $\mu$ on $\O_{K,\mathfrak l}$ which is supported on $\O_{K,\mathfrak l}^\times$, we attach a new measure $\beta=[\rho_{\varphi}^{-1}]\mu$ on $\O_{K,\mathfrak l}$, which is also supported on $\O_{K,\mathfrak l}^\times$, defined as follows; For all locally constant function $f$ on $\O_{K,\mathfrak l}$, we define $\int_{\O_{K,\mathfrak l}}f(x)d\beta(x):=\int_{\O_{K,\mathfrak l}}\rho_{\varphi}(x)^{-1}f(x)d\mu(x)$. Here, $\rho_{\varphi}^{-1}:\O_{K,\mathfrak l}\rightarrow \mathcal{W}$ is extended by zero.
Choose a finite unramified extension $F/\mathbf Q_p$ over which $R$ is defined and all points of $E'[\ell]$ are rational and which contains $\mu_\ell$ and the values of $\rho_\varphi$. Such an $F$ exists by Theorem~\ref{thm:rational-representative} and the definition of $R$ in \eqref{eq:lm-theta}. The measure $\mu$ takes values in $F(\mu_{\ell^\infty})$, since $E'[\mathfrak l^\infty]$ is rational over this field \cite[\S2]{Lam15}. The finite Fourier expansion of $\rho_\varphi^{-1}$ on $\O_K/\mathfrak l$ expresses the Fourier transform of $\beta$ as an $F$-linear combination of translates of $R$ by points of $E'[\mathfrak l]$. Thus $\beta$ is an elliptic function measure for $E'/F$.
Since $\mu$ is supported on the units, twisting $\beta$ by $\rho_\varphi$ recovers $\mu$. Applying \cite[Proposition~3.4]{Lam15} in both directions gives $v(\beta)=v(\mu)=v(R)=0$. Moreover, $\rho_\varphi$ is trivial on $\mathcal O_K^\times$, so
$\beta$ is invariant under these units. Lamplugh's theorem \cite[Theorem~5.1]{Lam15} therefore gives
\begin{equation}\label{eq:gamma-nonvanishing}
 \Gamma_\mu(\rho_{\varphi\psi}^{-1})
   =\Gamma_\beta(\rho_\psi^{-1})\in\mathcal W^\times
\end{equation}
for all but finitely many $\psi$.

The remaining interpolation factors are units for all but finitely
many $\psi$. Indeed, if a nonunit $c\in\mathcal O_K$ is prime to
$\ell$, its image in $\Gamma_{\mathfrak l}$ is nontrivial:
otherwise its image in $K_{\mathfrak l}^\times$ would be a root of
unity, forcing $c$ itself to be a root of unity.
Consequently, the order of $\psi((c))$ tends to infinity with the
conductor of $\psi$. Since reduction is injective on
$\ell$-power roots of unity, $\psi((c))$ can equal any fixed
nonzero residue value for only finitely many $\psi$.
Applying this observation to $c=N$ shows that
$B_N(\lambda\varphi\psi)$ is a unit if the conductor of $\psi$ is sufficiently large.
It also applies to generators of the primes in
\begin{equation}\label{eq:lm-tame-euler}
 L_{\mathfrak g}^{\mathrm{alg}}(\lambda\chi,0)
   =\prod_{\substack{\mathfrak r\mid\mathfrak g\\
                      \mathfrak r\nmid\mathfrak f_\lambda}}
       \bigl(1-(\lambda\chi)(\mathfrak r)\bigr)
      L^{\mathrm{alg}}(\lambda\chi,0).
\end{equation}
Now, take $\chi=\varphi\psi$ in
Proposition~\ref{prop:lm-interpolation}. Equations
\eqref{eq:gamma-nonvanishing} and \eqref{eq:lm-tame-euler}, together
with the unit interpolation factors, prove
\eqref{eq:main-nonvanishing}.

\section*{Acknowledgements}
I thank Antonio Lei for his feedback.
This work was supported by JSPS KAKENHI Grant Number JP25KJ1217.
I used OpenAI's ChatGPT (\url{https://chatgpt.com}, November 2026) to help phrase a few passages in the proofs.

\end{document}